\documentclass[final, sort&compress]{elsarticle} 
\usepackage[text={13.5cm, 22cm}, centering] {geometry} 

\usepackage{amsmath,amsthm,amsfonts,amssymb}
\usepackage[all,poly,knot]{xy}

\usepackage{longtable,tabularx,booktabs}
\usepackage[flushleft]{threeparttable}
\allowdisplaybreaks[4]
\usepackage{setspace} 
\usepackage{multirow}
\usepackage{pdflscape}  

\usepackage{xcolor}
\usepackage[colorlinks=true,allcolors=teal]{hyperref} 

\usepackage[nameinlink,capitalize]{cleveref} 
\crefformat{equation}{#2(#1)#3}
\crefmultiformat{equation}{(#2#1#3)}{ and~(#2#1#3)}{, (#2#1#3)}{ and~(#2#1#3)}

\newtheorem{theorem}{Theorem}
\newtheorem{corollary}{Corollary}
\newtheorem{lemma}{Lemma}

\newtheorem{example}{Example}

\newcommand{\Z}{\mathbb{Z}}
\newcommand{\F}{\mathbb{F}}
\newcommand{\N}{\mathrm{N}}

\begin{document}
\begin{frontmatter}
\title{On inverse of permutation polynomials of small degree \\over f{i}nite f{i}elds, II\tnoteref{t1}}

\author[QF]{Yanbin Zheng}
\ead{zhengyanbin16@126.com}

\author[GZ]{Yuyin Yu}
\ead{yuyuyin@163.com}

\address[QF]{School of Mathematical Sciences, Qufu Normal University, Qufu 273165, China}
\address[GZ]{School of Mathematics and Information Science, Guangzhou University, Guangzhou 510006, China} %

\begin{abstract}
  We investigate the permutation property of polynomials of the form $x^{r}(x^{s} -a)^{t}$,
  and give the expressions of their inverses.
  In particular, explicit expressions of inverses of permutation polynomials
  $x(x^3 -a)^2$ and $x(x^2 -a)^3$ on~$\mathbb{F}_{7^n}$ are presented.
  Then, using some known results, we obtain the inverses of all permutation polynomials
  of degree $6, 7, 8$ over finite fields.
\end{abstract}

\end{frontmatter}
\section{Introduction}
For $q$ a prime power, let $\F_{q}$ denote the finite field with $q$ elements,
$\F_{q}^{*} = \F_{q} \setminus \{0\}$, and $\F_{q}[x]$ the ring of polynomials over~$\F_{q}$.
A polynomial $f \in \F_{q}[x]$ is called a \emph{permutation polynomial} (PP) of~$\F_{q}$
if it induces a bijection from~$\F_{q}$ to itself. For any PP $f$ of~$\F_{q}$,
there exists a polynomial $f^{-1} \in \F_{q}[x]$ such that $f^{-1}(f(c)) =c$
for each $c \in \F_{q}$ or equivalently $f^{-1}(f(x)) \equiv x \pmod{x^{q} -x}$,
and the polynomial $f^{-1}$ is unique in the sense of reduction modulo $x^q -x$.
Hence $f^{-1}$ is defined as the \emph{composition inverse} of~$f$,
and we simply call it the \emph{inverse} of $f$ on~$\F_{q}$.

A polynomial over $\F_{q}$ is called an \emph{exceptional polynomial} over~$\F_{q}$
if it is a PP of~$\F_{q^n}$ for infinitely many positive integers~$n$.
Two polynomials $f$ and $g$ over $\F_{q}$ are called
\emph{affine equivalence}
if $g(x) = \alpha f(\beta x + \gamma) + \delta$ for some
$\alpha, \beta \in \F_{q}^{*}$ and $\gamma, \delta \in \F_{q}$.
Affine equivalent $f$ and $g$ share the same degree,
and $f$ is a PP of $\F_{q}$ if and only if so is $g$.

The classification of PPs of finite fields has a long history.
In~1896, Dickson~\cite{Dickson1896} obtained all \emph{normalized} PPs of degree
$\leq 5$ of~$\F_{q}$ for all $q$, and classified all PPs of degree~$6$ of $\F_{q}$ for odd~$q$.
In~2010, a complete classification of all PPs of degree~$6$ or~$7$ of $\F_{2^n}$
was settled in~\cite{JLi10}, up to affine equivalence and a special transformation.
However, each class of resulting PPs is invariant
under the special transformation~\cite{Fan7}.
For a verification of the classification of
normalized PPs of degree $6$ of $\F_{q}$ for all $q$, see~\cite{PPs6}.
More recently, under affine equivalence, Fan \cite{Fan7,Fan8even,Fan8odd}
gives a complete classification of all PPs of degree~$7$ of~$\F_{q}$
for odd~$q$ and degree~$8$ of~$\F_{q}$ for all~$q$.
All such PPs of degree~$\leq 8$ can be divided into two classes:
exceptional and non-exceptional.
According to the results in the above literature,
a non-exceptional PP of degree~$\leq 8$ over $\F_{q}$ exists only if $q < 64$.

The inverses of all normalized PPs of degree $\leq 5$ were listed in~\cite{Zheng7deg}.
In this paper, we consider the inverses of all PPs of degree $6, 7, 8$.
For non-exceptional PPs of degree $6, 7, 8$ of $\F_{q}$, since $q < 64$ is very small,
one can use the Lagrange interpolation formula to compute their inverses, i.e.,
\begin{equation}\label{Lag-1}
f^{-1}(x) = \sum_{c \in \F_{q}}c \big(1-(x-f(c))^{q-1}\big).
\end{equation}
But for exceptional PPs of degree $6, 7, 8$ of $\F_{q^n}$ with infinitely many~$n$,
finding explicit expressions of their inverses on $\F_{q^n}$ is not an easy problem.

There are several papers on the inverses of some classes of PPs,
see for example \cite{Wu-L-bi,Wu-L,Reis18} for linearized PPs,
\cite{KLi18-1,MR-1,Wang-1,ZhengPW2} for PPs of the form $x^rh(x^{(q-1)/d})$,
\cite{Involutions,NiuLQW19,ZhengYL+2019} for involutions over $\F_{2^n}$,
\cite{ZhengPW2,Wang-cyc2} for generalized cyclotomic mapping PPs,
\cite{ZhengPW1,ZhengPW2,ZhengRDP} for more general piecewise PPs,
\cite{TW-1,TW17,NiuLQW20} for PPs constructed by the AGW criterion.
The results in \cite{TW-1,TW17,NiuLQW20} contain some concrete classes
such as bilinear PPs \cite{Coulter-1,Wu-bil},
linearized PPs of the form $L(x)+K(x)$\cite{Reis18},
and PPs of the form $x+ \gamma f(x)$ \cite{Kyur11}.
For a brief summary of the results concerning the inverses of PPs,
we refer the reader to~\cite{Zheng7deg} and the references therein.

In this paper, we study the PPs of the form $f(x) = x^{r}(x^{s} - a)^{t}$
on $\F_{q^n}$, where~$a \in \F_{q^n}^{*}$, $st=q^m-1$,
$r \equiv 1 \pmod{\ell}$ and $\ell=(q^{n}-1)/(s, q^n -1)$.
According to the Akbary-Ghioca-Wang (AGW) criterion \cite{AGW},
$f$ is a PP of $\F_{q^n}$ if and only if $(r, s, q^n-1) = 1$ and another
polynomial $g(x) := x(x-a)^{st}$ permutes the subset $(\F_{q}^{*})^s$.
By solving the equation $g(x) = c$ for any $c \in (\F_{q}^{*})^s$,
we find the inverse $g^{-1}$ of $g$ on $(\F_{q}^{*})^s$, and prove that
$g$ permutes $(\F_{q}^{*})^s$ if and only if $a^\ell \neq 1$.
Substituting $g^{-1}$ into a slightly modified version of
a result in \cite{NiuLQW20} concerning the inverse of
more general PP $x^r h(x^s)$,
we obtain an expression of the inverse of $f$ on $\F_{q^n}$.

By considering special cases such as
$m=n$, $q^{(m,n)} -1 \mid s$, and $t=2$ or $3$,
we get some new classes of PPs and their inverses.
In particular, explicit expressions of inverses of exceptional polynomials
$x(x^3 -a)^2$ and $x(x^2 -a)^3$ on $\F_{7^n}$ are given.
Then, based on the known formulae for the inverses of Dickson PPs and linearized PPs,
we find the inverses of all exceptional polynomials of degree $6, 7, 8$; see~\cref{T1}.

In summary, under affine equivalence,
\cref{T1} and~\cite[Tabel~I]{Zheng7deg} list
the inverses of all PPs of degree $\leq 8$ over all finite fields,
except for the inverses of non-exceptional PPs of degree $6, 7, 8$
which can be obtained by~\cref{Lag-1}.

Some notations of this paper are as follows.
The sets of integers and positive integers are denoted by
$\Z$ and $\mathbb{N}$ respectively.
The greatest common divisor of two integers $m$ and $n$ is written as $(m, n)$.
For $a \in \F_{q^n}$ and $d \mid n$, the norm of $a$ over $\F_{q^d}$
is defined by $\N_{q^n/q^d}(a) = a^{(q^n-1)/(q^d-1)}$.

\begin{landscape}
\begin{table*} [t]
\caption{All exceptional polynomials of degree $6, 7, 8$ and their inverses}
\centering
\begin{threeparttable}
\begin{spacing}{1.7}
\begin{tabular}{l l l l}\toprule[1pt] \label{T1}
  Exceptional polynomials over $\F_{q}$
         & Inverses & $q$  & Reference \\\midrule
  $x^6$  & $x^{(5q-4)/6}$
         & $q = 2^{n}$, odd $n \geq 3$
         & \cref{lem:xm}\\
  $x^7$  & $x^{(kq-k+1)/7}$ with $k \equiv (1-q)^{5} \pmod{7}$
         & $q \not\equiv 1 \pmod 7$
         & \cref{lem:xm}\\
  $x^{7} - ax$   ($a$ not a sixth power)
        & $a^{\frac{q-1}{6}}(1-a^{\frac{q-1}{6}})^{-1}
          \sum_{i=0}^{n-1}a^{-\frac{7^{i+1} -1}{6}}x^{7^{i}}$
        & $q = 7^{n}$, $n \geq 2$
        & {\cite{CoulterH04,Wu-L-bi}} \\
  $x^{7} -2 a x^{4} +a^2 x$ \;($a$ not a cube)
        & $
            2x\big(2 a^{\frac{q-1}{6}} x^{\frac{q-1}{2}} + a^{\frac{q-1}{3}} + 1\big)
          \big(\sum_{i=0}^{n-1} a^{-\frac{7^{i+1}-1}{6}}x^{\frac{7^{i}-1}{2}}\big)^{2}$
        &$q = 7^n$, $n \geq 2$
        &\cref{cor:r=1+t=2+x7}\\
  $x^{7} -3 a x^{5} +3a^2 x^3 -a^3 x$ \;($a$ not a square)
        & $x\big( 3 (a x^4)^{\frac{q-1}{6}} -3 (a x)^{\frac{q-1}{3}}  -2 \big)
        \big(\sum_{i=0}^{n-1} a^{-\frac{7^{i+1}-1}{6}}x^{\frac{7^{i}-1}{3}}\big)^{3}$	
        &$q = 7^n$, $n \geq 2$
        & \cref{cor:r=1+t=3+x7}\\
  $x^7 - 7a x^5 + 14 a^2 x^3 - 7 a^3 x$ \;($a \neq 0$)
        & $\sum_{i=0}^{\lfloor m/2\rfloor}
             \frac{m}{m-i} \binom{m-i}{i}(-a^7)^i x^{m-2i}$
        &$q \equiv \pm 2, \pm 3 \pmod{7}$
        &{\cite{KLi18-1}}\\
        & where $m = (kq^{2} - k + 1)/7$ with $k \equiv (1 - q^{2})^5 \pmod{7}$
        &&\cref{lem:xm}\\
        & and $\lfloor m/2\rfloor$ denotes the largest integer $\leq m/2$.  &&\\
  $x^{8} + a_{2} x^{4} + a_{1} x^{2}  + a_{0} x$
        & $(\det(D_{L}))^{-1}\sum_{i=0}^{n-1} \bar{a}_i x^{q^i}$
        &$q = 2^n$, $n \geq 4$
        &{\cite{Wu-L}}\\
  (if its only root in $\F_{2^n}$ is~$0$)
        & where $D_{L}$ and $\bar{a}_i$ are defined as in \cref{lem:L}.  &&\\
  \bottomrule[1pt]
  \end{tabular}
\end{spacing}
\vspace{-1ex}
  \begin{tablenotes}
    \item[$\dagger$] This list is complete up to affine transformations:
                      $g(x) = \alpha f(\beta x + \gamma) + \delta$
                      with $\alpha,\beta,\gamma,\delta \in \F_{q}$ and $\alpha\beta \neq 0$.
    \item[$\ddagger$] All non-exceptional permutation polynomials of degree $6, 7, 8$ are listed in~\cite{JLi10,PPs6,Fan7,Fan8even,Fan8odd}, and
        all of them are over small fields $\F_{q}$ with $q \leq 64$. Hence
        their inverses can be obtained by the Lagrange interpolation formula.
  \end{tablenotes}

\end{threeparttable}
\end{table*}
\end{landscape}

\section{Permutation polynomials of small degree and their inverses}

Since the inverses of all normalized PPs of degree $\leq 5$ were listed in~\cite{Zheng7deg},
we only consider the PPs of degree $6, 7, 8$
which can be divided into two classes: exceptional and non-exceptional.
Combining the results in~\cite{Bartoli8even,JLi10,PPs6,Fan7,Fan8even,Fan8odd}
gives the following theorems.

\begin{theorem}\label{thm:nEPs}
A non-exceptional PP of degree $m \in \{6, 7, 8\}$ over $\F_{q}$
exists if and only if one of the following conditions holds:
\begin{enumerate}[\upshape(i)]
  \item  $m = 6$ and $q \in \{8, 9, 11, 16, 27, 32\};$
  \item  $m = 7$ and $q \in \{9, 11, 13, 16, 17, 19, 23, 25, 27, 31, 49\};$
  \item  $m = 8$ and $q \in \{11, 13, 16, 19, 23, 27, 29, 31, 32, 64\}.$
\end{enumerate}
Under affine equivalence, all such PPs are explicitly listed in
\textnormal{\cite{JLi10,PPs6,Fan7,Fan8odd,Fan8even}}.
\end{theorem}

\begin{theorem}\label{thm:EPs}
Each exceptional polynomial of degree $6, 7, 8$
is affine equivalent to one of the following:
\begin{enumerate}[\upshape(i)]
  \item $x^6$ over $\F_{2^n}$\! with odd $n \geq 3;$
  \item $x^7$ over $\F_{q}$\! with $q \not\equiv 1 \pmod{7};$
  \item $x^{7} - ax$ over $\F_{7^n}$\! with $a$ not a sixth power in $\F_{7^n}$, i.e.,
        $a \in \F_{7^n}^{*}$\! such that $a^{(7^n-1)/6} \neq 1;$
  \item $x(x^{3} - a)^2$ over $\F_{7^n}$\! with $a$ not a cube in $\F_{7^n}$, i.e.,
        $a \in \F_{7^n}^{*}$\! such that $a^{(7^n-1)/3} \neq 1;$
  \item $x(x^{2} - a)^3$ over $\F_{7^n}$\! with $a$ not a square in $\F_{7^n}$, i.e.,
        $a \in \F_{7^n}^{*}$\! such that $a^{(7^n-1)/2} \neq 1;$
  \item $x^7 - 7a x^5 + 14 a^2 x^3 - 7 a^3 x$
         with $a \in \F_{q}^{*}$\! and $q \equiv \pm 2, \pm 3 \pmod{7};$
  \item $x^{8} + a_{2} x^{4} + a_{1} x^{2}  + a_{0} x$ over $\F_{2^n}$\!
        if its only root in $\F_{2^n}$\! is~$0$.
\end{enumerate}
\end{theorem}

The inverses of PPs in~\cref{thm:nEPs} can be obtained directly
by the Lagrange interpolation formula~\cref{Lag-1} due to $q < 64$.
To obtain the inverses of PPs in~\cref{thm:EPs},
we need the following results. 

\begin{lemma}\label{lem:xm}
Let $m, \ell \in \mathbb{N}$ and $(m, \ell)=1$.
Then an inverse of~$m$ modulo~$\ell$ is $(k \ell + 1)/m$, where
$k \equiv -\ell^{\phi(m)-1} \pmod{m}$
and~$\phi$ is Euler's totient function.
\end{lemma}
\begin{proof}
Clearly, $k \ell + 1 \equiv 1 - \ell^{\phi(m)} \equiv 0 \pmod{m}$
and $m(k \ell + 1)/m - k \ell =1$.
Therefore, $(k \ell + 1)/m$ is an inverse of~$m$ modulo~$\ell$.
\end{proof}

\begin{lemma}[{\cite{CoulterH04,Wu-L-bi}}]\label{lem:L-bi}
  Let $L(x) =x^{q^m} -ax$, where $a \in \F_{q^n}^{*}$ and $m, n \in \mathbb{N}$.
  Then $L$ is a PP of $\F_{q^n}$\! if and only if $\N_{q^n/q^d}(a) \neq 1$,
  where $d =(m,n)$. In this case, its inverse on $\F_{q^n}$ is
  \[
    L^{-1}(x)=\frac{\N_{q^n/q^d}(a)}{1-\N_{q^n/q^d}(a)}
    \sum_{i=1}^{n/d}a^{-\frac{q^{im} -1}{q^m-1}}x^{q^{(i-1)m}}.
  \]
\end{lemma}

\begin{lemma}\label{lem:L}
Let $L(x) = \sum_{i=0}^{n-1} a_i x^{q^i} \in \F_{q^n}[x]$.
Then $L$ is a PP of $\F_{q^n}$\! if and only if
\[
D_{L}  :=
\left(
  \begin{array}{cccc}
 a_{0}       & a_{1}     &  \cdots  & a_{n-1}      \\
 a_{n-1}^{q} & a_{0}^{q} &  \cdots  & a_{n-2}^{q}  \\
 \vdots      & \vdots    &  \vdots  & \vdots       \\
 a_{1}^{q^{n-1}}  & a_{2}^{q^{n-1}} &  \cdots  & a_{0}^{q^{n-1}}
  \end{array}
\right)
\]
is nonsingular~\textnormal{\cite[Page 362]{FF}}.
In this case, its inverse was given in~\textnormal{\cite{Wu-L}} by
\[
L^{-1}(x) = \frac{1}{\det(D_{L})}
\sum_{i=0}^{n-1} \bar{a}_i x^{q^i},
\]
where $\bar{a}_i$ is the $(i,0)$-th cofactor of $D_{L}$, i.e.,
$
\det(D_{L}) = a_0 \bar{a}_0 + \sum_{i=1}^{n-1} a_{n-i}^{q^i} \bar{a}_i.
$
\end{lemma}

The Dickson polynomial $D_{n}(x,a)$ of degree $n$
with parameter $a \in \F_{q}$ is given as
\[
D_{n}(x,a) = \sum_{i=0}^{\lfloor n/2\rfloor}
             \frac{n}{n-i} \binom{n-i}{i}(-a)^i x^{n-2i},
\]
where $\lfloor n/2\rfloor$ denotes the largest integer $\leq n/2$.
For $a \in \F_{q}^{*}$, $D_{n}(x,a)$ is a PP of $\F_q$ if and only if $\gcd(n ,q^2 -1)=1$.
Its inverse is determined in~\cite{KLi18-1} by the next lemma.

\begin{lemma}[{\cite[Lemma 4.8]{KLi18-1}}]\label{lem:dickson-1}
  Let $a \in \F_{q}^{*}$ and $m, n \in \mathbb{N}$ be such that $mn \equiv 1 \pmod{q^2-1}$.
  Then the inverse of $D_{n}(x,a)$ on $\F_q$ is $D_{m}(x, a^n)$.
\end{lemma}

The PP $x^7 - 7a x^5 + 14 a^2 x^3 - 7 a^3 x$ in \cref{thm:EPs} is 
the Dickson polynomial $D_{7}(x, a)$ and,
by \cref{lem:dickson-1}, its inverse is $D_{m}(x, a^7)$.
By \cref{lem:xm}, an inverse $m$ of $7$ modulo $q^{2}-1$ can be written as
$m = (kq^{2} - k + 1)/7$, where $k \equiv (1 - q^{2})^5 \pmod{7}$.

The inverses of PPs $x^6$, $x^7$, $x^7 -ax$ and $x^8 + \sum_{i=0}^{2}a_{i}x^{i}$
in \cref{thm:EPs} can be obtained directly by \cref{lem:xm,lem:L-bi,lem:L}.
The inverses of PPs $x(x^{3} - a)^2$ and $x(x^{2} - a)^3$ on $\F_{7^n}$
are given by \cref{cor:r=1+t=2+x7,cor:r=1+t=3+x7} in \cref{sec:PPs34}
respectively.

In short, the inverses of all the PPs in \cref{thm:EPs} are known now.
For convenience, all the PPs in \cref{thm:EPs} and their inverses
are listed in \cref{T1}.

\section{Large class of PPs and its inverse}

In this section we slightly modify a known expression of the inverse of
PP $x^rh(x^{s})$. Moreover, we investigate the permutation properties of
$x(x-a)^{st}$ on $(\F_{q^n}^{*})^s$ and $x^{r}(x^{s} -a)^{t}$ on $\F_{q^n}$,
and obtain their inverses. The following lemma will be needed.

\begin{lemma}[{\cite[Proposition 3.1]{AGW}}]\label{lem:xrh(xs)}
  Let $f(x) = x^rh(x^{s})$, where $h \in \F_{q}[x]$ and $r, s \in \mathbb{N}$.
  Then $f$ is a PP of $\F_{q}$ if and only if $(r, s, q-1) = 1$
  and $g(x) := x^r h(x)^{s}$ permutes $(\F_{q}^{*})^s$.
\end{lemma}

The problem that when special classes of $g$ permuting $(\F_{q}^{*})^s$
has been extensively studied, see for example
\cite{Bartoli20-tri,CaoHMX20,HouTZ20-tri,KLi18,KLiQLC20,NLiH19,TuLZ19,ZhaHZ19}.
For a recent survey of this problem,
we refer the reader to \cite{WangIndex19} and the references therein.

The inverse of $f$ in terms of roots of unity
over $\F_{q}$ was given in \cite{Wang-1,ZhengPW2}.
Let $g_1 \in \F_{q}[x]$ be such that
$x^{k_1} \circ g_1 = g \circ x^{k_1}$, where $k_1 = s/(s,q-1)$.
The inverse of $f$ in terms of the inverse of $g_1$ on $(\F_{q}^{*})^{s}$
was given in~\cite{KLi18-1} when $(r, q-1)=1$,
and given in~\cite{NiuLQW20} for all $r \in \mathbb{N}$.
Using the method in~\cite{NiuLQW20}, we obtain the following equivalent version of
the inverse of $f$, which is expressed in terms of
the inverse of $g$ on $(\F_{q}^{*})^s$.

\begin{theorem}\label{thm:xrh(xs)inv2}
  Let $f(x) = x^rh(x^{s})$, where $h \in \F_{q}[x]$ and $r, s \in \mathbb{N}$.
  Let $\bar{s} = (s, q-1)$ and $k$ be an inverse of
  $s/\bar{s}$ modulo~$(q-1)/\bar{s}$.  If~$f$ is a PP of $\F_{q}$ and
  $g^{-1}$ is the inverse of $g(x) := x^r h(x)^{s}$ on $(\F_{q}^{*})^s$.
  Then the inverse of~$f$ on $\F_{q}$ is
  \[
   f^{-1}(x) = x^{b} \big(h(g^{-1}(x^s))\big)^{-b} \big(g^{-1}(x^s)\big)^{ck},
  \]
  where $b, c \in \Z$ satisfy $br + c\bar{s}= 1$.
\end{theorem}
\begin{proof}
  We prove it by using the method in {\cite[Theorem~3.2]{NiuLQW20}}.
  Since the theorem holds for $x=0$, we only consider $x \in \F_{q}^{*}$.
  Let $\phi(x) = (x^{r}, x^{s})$ and $\psi(y, z) = (y^{r}h(z)^{r}, g(z))$.
  It is easy to verify that $\psi \circ \phi = \phi \circ f$, i.e.,
  the following diagram is commutative:
\begin{equation*}
\xymatrix{
   \F_{q}^{*} \ar[r]^{f}\ar[d]_{\phi}   & \F_{q}^{*} \ar[d]^{\phi} \\
   \phi(\F_{q}^{*}) \ar[r]^{\psi}       & \phi(\F_{q}^{*})
.}
\end{equation*}

If $f$ is a PP, then, by \cref{lem:xrh(xs)}, $(r, \bar{s}) = 1$ and
$g$ permutes $(\F_{q}^{*})^s$. Assume $br + c\bar{s}= 1$
and $k (s/\bar{s}) + v (q-1)/\bar{s} =1$  for some $b, c, k, v \in \Z$. Then
\[
  c\bar{s} = c\bar{s}(k (s/\bar{s}) + v (q-1)/\bar{s})
  = cks + cv(q-1),
\]
and so, for any $x \in \F_{q}^{*}$,
\begin{equation}\label{eq:xrbxsck=x}
  (x^r)^b (x^s)^{ck}
  = x^{br + cks} x^{cv(q-1)}
  = x^{br + cks + cv(q-1)}
  = x^{br + c\bar{s}}
  = x.
\end{equation}
Hence $\phi$ is bijective and $\phi^{-1}(y, z) = y^b z^{ck}$.

Since $f$ is a PP of $\F_{q}$ and $\phi$ is bijective, $\psi$ is also bijective.
For $(y, z) \in \phi(\F_{q}^{*})$, assume that $\psi(y, z) = (\alpha, \beta)$
for some $(\alpha, \beta) \in \phi(\F_{q}^{*})$, i.e.,
\[
y^r h(z)^r = \alpha \quad \text{and} \quad
g(z) = z^r h(z)^s = \beta.
\]
Because $z, \beta \in (\F_{q}^{*})^s$ and $g$ permutes $(\F_{q}^{*})^s$,
we obtain $z = g^{-1}(\beta)$.
Denote $\alpha = x_{0}^{r}$ and $\beta = x_{0}^{s}$
for some $x_{0} \in \F_{q}^{*}$. Then $\alpha^b \beta^{ck} = x_{0}$ by~\cref{eq:xrbxsck=x}, and so
$
  \big( \alpha^b \beta^{ck} h(z)^{-1}\big)^r h(z)^r = \alpha
$.
Since $\psi$ is bijective, we have $y = \alpha^b \beta^{ck} h(z)^{-1}$.
Hence,
\begin{equation*}\label{eq:psi-inv}
  \psi^{-1}(\alpha, \beta) = (y, z)
  = (\alpha^b \beta^{ck} h(g^{-1}(\beta))^{-1}, g^{-1}(\beta)).
\end{equation*}
Substituting $\phi, \phi^{-1}, \psi^{-1}$  into
$f^{-1} = \phi^{-1} \circ \psi^{-1}  \circ  \phi$
gives the desire result.
\end{proof}

The key step in \cref{thm:xrh(xs)inv2} is to find
the inverse of $g$ on $(\F_{q}^{*})^s$,
which is possible to be done for special classes of $g$,
as for instance in the following result.

\begin{lemma} \label{lem:Subset}
  Let $a \in \F_{q^n}^{*}$\! and $s \mid q^m-1$.
  Then $g(x) = x(x-a)^{q^m-1}$ permutes $(\F_{q^n}^{*})^s$ if and only if
  $a^{\ell} \neq 1$, where $\ell = (q^{n}-1)/(s, q^n -1)$.
  In this case, its inverse on $(\F_{q^n}^{*})^s$ is
  \[
   g^{-1}(x) = \Big((a^{-1}x)^{\frac{q^n-1}{q^{d}-1}} -1\Big)
      \bigg(\sum_{i=1}^{n/d}a^{-\frac{ q^{i m}-1 }{q^m-1}}
       x^{\frac{q^{(i-1)m}-1}{q^m-1}}\bigg)^{-1} + a,
  \]
  where $d = (m, n)$.
\end{lemma}
\begin{proof}
Since the multiplicative group of $\F_{q^n}$ is cyclic,
we can verify $(\F_{q^n}^{*})^s = (\F_{q^n}^{*})^{(s, q^n -1)}$.
Thus $a \in (\F_{q^n}^{*})^s$ if and only if $a^{\ell} =1$.
If $a^{\ell} =1$, then $a \in (\F_{q^n}^{*})^s$, and so
$g$ has root in $(\F_{q^n}^{*})^s$. Hence $g$ does not permute $(\F_{q^n}^{*})^s$.
Next we only consider $a^{\ell} \neq 1$ and $x \in (\F_{q^n}^{*})^s$.

Since $g$ introduces a mapping from $(\F_{q^n}^{*})^s$ to itself,
we need only show that, for any $y \in (\F_{q^n}^{*})^s$,
the equation $g(x) = y$ has exactly one solution
$x$ and $x = g^{-1}(y)$.

Since $a^{\ell} \neq 1$, we have $x - a \neq 0$ for any $x \in (\F_{q^n}^{*})^s$.
Let $z = (x - a)^{-1}$. Then $x - a = z^{-1}$ and $x = z^{-1} + a$.
Substituting them into $g(x) =y$ yields
  \begin{equation}\label{eq:L(z)}
    (z^{-1} + a) z^{1-q^m} = y,  \quad \text{i.e.,} \quad
    z^{q^m} -(a/y) z = 1/y.
  \end{equation}
  Recall that $d =(m, n)$ and $\ell = (q^n-1)/(s, q^n-1)$. Let $q^m-1 = st$. Then
  \[\begin{split}
    q^d -1
    &= q^{(m,n)}-1 = (q^m-1,q^n-1) = (st,q^n-1) \\
    &= (s,q^n-1)(t,(q^n-1)/(s,q^n-1)) \\
    &= (s,q^n-1) (t, \ell),
  \end{split}\]
  and so $\frac{q^n-1}{q^d -1}(t, \ell)= \ell$.
  Since $a^{\ell} \neq 1$ and $y \in (\F_{q^n}^{*})^s$, we have
  \begin{equation*}
    \big(\N_{q^n/q^d}(a/y)\big)^{(t,\ell)}
    = (a/y)^{\frac{q^n-1}{q^d -1}(t, \ell)}
    = (a/y)^{\ell} = a^{\ell}  \neq 1.
  \end{equation*}
  Therefore,
  \begin{equation}\label{eq:Na-neq-Nys}
    \N_{q^n/q^d}(a/y) \neq 1.
  \end{equation}
  It follows from \cref{lem:L-bi} and~\cref{eq:L(z)} that
  \begin{equation}\label{eq:z=f(y)}
  \begin{split}
  z &= \frac{\N_{q^n/q^d}(a/y)}{1-\N_{q^n/q^d}(a/y)}
       \sum_{i=1}^{n/d} (a/y)^{-\frac{ q^{i m}-1 }{q^m-1} }
       (1/y)^{q^{(i-1)m}} \\
    &= \big(\N_{q^n/q^d}(y/a) - 1 \big)^{-1}
        \sum_{i=1}^{n/d} a^{-\frac{ q^{i m}-1 }{q^m-1}}
        y^{\frac{q^{(i-1)m}-1}{q^m-1}}.
  \end{split}
  \end{equation}
  Substituting~\cref{eq:z=f(y)} into $x = z^{-1} + a$ gives $x = g^{-1}(y)$.
\end{proof}

After these preparations, we can now give the main theorem.

\begin{theorem}\label{thm:rst}
Let $f(x)= x^r(x^{s} -a)^{t}$, where $a \in \F_{q^n}^{*}$\!
and $r, s, t \in \mathbb{N}$ are such that $st = q^m - 1$ and
$r \equiv 1 \pmod{\ell}$ with $\ell = (q^{n}-1)/(s, q^n -1)$.
Then $f$ is a PP of $\F_{q^n}$ if and only if
$(r, q^{n}-1)=1$ and $a^{\ell} \neq 1$.
In this case, the inverse of $f$ on $\F_{q^n}$ is
\begin{equation}\label{eq:rst-inv}
  f^{-1}(x) = x^{u}  (G(x) H(x))^{tu},
\end{equation}
where $u$ is an inverse of\, $r$  modulo $q^{n}-1$,
\begin{align}
    G(x) &= \big( (a^{-1}x^{s})^{\frac{q^n-1}{q^d-1}} -1 \big)^{-1}, \label{eq:gx}\\
    H(x) &= \sum_{i=1}^{n/d}
    a^{-\frac{ q^{i m}-1 }{q^m-1}}
    x^{\frac{q^{(i-1)m}-1}{t}},
    ~\text{and $d = (m, n)$.}  \label{eq:hx}
\end{align}
\end{theorem}
\begin{proof}
  From \cref{lem:xrh(xs)}, $f$ is a PP if and only if
  $(r, s, q^{n}-1)=1$ and $g(x) := x^r (x-a)^{st}$ permutes $(\F_{q}^{*})^s$.
  Since $r \equiv 1 \pmod{\ell}$, we have $x^r = x$
  and so $g(x) = x (x-a)^{st}$ for any $x \in (\F_{q}^{*})^s$.
  By \cref{lem:Subset},
  $g$ permutes $(\F_{q}^{*})^s$ if and only if $a^{\ell} \neq 1$.
  It follows from $r \equiv 1 \pmod{\ell}$ that $(r, \ell) = 1$,
  and so $(r, s, q^{n}-1)=1$ if and only if $(r, q^{n}-1)=1$.

  Let $ru + (q^{n}-1)v = 1$ for some $u, v \in \Z$.
  Then $ru + (\ell v) \bar{s} = 1$, where $\bar{s} = (s, q^{n}-1)$.
  For any $x \in \F_{q^n}^{*}$, we have $g^{-1}(x^s) \in (\F_{q}^{*})^s$,
  and so  $(g^{-1}(x^s))^{\ell} =1$.
  Substituting $b := u$, $c := \ell v$ and $g^{-1}$ in \cref{lem:Subset}
  into \cref{thm:xrh(xs)inv2} gives the expression of $f^{-1}$.
\end{proof}
The following are two examples of \cref{thm:rst}
where $r^2 \equiv 1 \pmod{q^n-1}$.

\begin{example}\label{Exam:rst-q=2}
Let $f(x)= x^r(x^{3} + a)^{5}$, where $a \in \F_{2^8}^{*}$\! and $r \equiv 1 \pmod{85}$.
Then $f$ is a PP of $\F_{2^8}$ if and only if $(r, 3) = 1$ and $a^{85} \neq 1$.
In this case, its inverse on $\F_{2^8}$ is
\[
f^{-1}(x) = x^{r}(x^{51} + a^{17})^{-5r}(x^{3} + a^{16})^{5r}.
\]
\end{example}
\begin{example}\label{Exam:rst-q=31}
Let $f(x)= x^r(x^{16} -a)^{5}$,
where $a \in \F_{3^6}^{*}$\! and $r \equiv 1 \pmod{91}$.
Then $f$ is a PP of $\F_{3^6}$ if and only if\, $(r, 8) = 1$ and $a^{91} \neq 1$.
In this case, its inverse on $\F_{3^6}$ is
\[
f^{-1}(x) = (1-a^{91})^{-5r} x^{r}
\big(a^{90} + a^{9}x^{16} + x^{584}\big)^{5r}.
\]
\end{example}

\section{Simplified versions of the main theorem}

In this section we aim to simplify the expressions of $G$ and $H$ in \cref{thm:rst}.
On the one hand, we consider small $n/(m, n)$ which is the number of the terms of $H$.
On the other hand, we study the cases $q^{d} - 1 \mid s$ and $q^{d}-1 \mid t$,
in which $G$ and $G^t$ are reduced to constants respectively.

\subsection{\bf The case $H$ is a constant}
Applying \cref{thm:rst} to $m = n$, we obtain that $G(x) = a(x^{s} -a)^{-1}$ and
$H(x) = a^{-1}$. Hence we arrive at the following result.

\begin{corollary} \label{cor:m=n}
Let $f(x)= x^r(x^{s} -a)^{t}$, where $a \in \F_{q}^{*}$,
$st = q - 1$ and $r \equiv 1 \pmod{t}$.
Then~$f$ is a PP of $\F_{q}$ if and only if $(r, q-1)=1$ and $a^{t} \neq 1$.
If $f$ is a PP of $\F_{q}$ and $u$ is an inverse of\, $r$  modulo~$q-1$,
then its inverse on $\F_{q}$ is
\begin{equation*}
  f^{-1}(x) =  x^{u} (x^{s} -a)^{-tu}.
\end{equation*}
\end{corollary}


The binomial $x^{s} -a$ in \cref{cor:m=n} can be generalized
to $H(x^s)$ which has no nonzero root in $\F_{q}$;
see the next theorem.


\begin{theorem}\label{thm:gen-m=n}
Let $f(x)= x^{r} (h(x^s))^{t}$, where $h \in \F_{q}[x]$ and $st = q-1$.
Then $f$ is a PP of $\F_{q}$ if and only if
$(r, q-1)=1$ and $h(x^s)$ has no nonzero root in $\F_{q}$
(see \cite[Corollary~3.2]{AW07}).
If $f$ is a PP of $\F_{q}$ and $r \equiv 1 \pmod{t}$,
then its inverse on $\F_{q}$ is
\[
f^{-1}(x) = x^{u} (h(x^s))^{-tu},
\]
where $u$ is an inverse of\, $r$  modulo~$q-1$.
\end{theorem}
\begin{proof}
The permutation part is \cite[Corollary~3.2]{AW07}.
Let $r = kt + 1$ for some $k \in \Z$.
Then
\begin{gather*}
  rs  = (kt + 1)s  = kst + s \equiv s \pmod{q-1}.
\end{gather*}
  Now it is easy to verify that
  $f^{-1}(f(e)) = e$ for any $e \in \F_{q}$.
  This completes the proof.
\end{proof}

%

\subsection{\bf The case $G$ or $G^t$ is a constant}
\vspace{3pt}
If $f$ in \cref{thm:rst} is a~PP, then~\cref{eq:Na-neq-Nys} holds,
and so $G(c) \in \F_{q^d}^{*}$ for any $c \in \F_{q^n}$.
Hence $G(x)^t = 1$ when $q^{d}-1 \mid t$. Moreover, if $q^{d} - 1 \mid s$,
then $\N_{q^n/q^d}(c^{s}) = 1$ for any $c \in \F_{q^n}^{*}$.
Thus $G$ is reduced to a constant. The argument above gives the following theorem.

\begin{theorem}\label{thm:qd-1-s-t}
\!With the same notation and hypothese as in \textnormal{\cref{thm:rst}},
let $f$ be a~PP~of\,~$\F_{q^n}$.
\begin{enumerate}[\upshape(i)]
  \item If $q^{d}-1 \mid t$, then 
$f^{-1}(x) = x^{u} (H(x))^{tu}$.
  \item If $q^{d}-1 \mid s$, then 
$f^{-1}(x) = x^{u} (AH(x))^{tu}$, where
$A = \big(a^{-\frac{q^n-1}{q^d-1}} - 1\big)^{-1}$.
\end{enumerate}
\end{theorem}

Recall that $d = (m, n)$ and $q^m -1 = st$.
The conditions $q^d -1 \mid s$ or $q^d -1 \mid  t$ in \cref{thm:qd-1-s-t}
are easy to satisfied, because
\begin{equation}\label{eq:st=qd-1}
  st = q^m -1 = (q^d -1)(1 + q^{d} + q^{2d} + \cdots + q^{m-d}).
\end{equation}
For instance, taking $t=1$ and $n=2m$ leads to
$q^d -1 = q^m -1 = s$ and $n/d=2$. Hence $H(x)=a^{-1} + a^{-(q^m+1)}x^{q^m-1}$.
Then substituting $q$ for $q^m$ yields the following result.

\begin{corollary} \label{cor:t=1&n=2m}
Let $f(x)= x^r(x^{q-1} -a)$, where $a \in \F_{q^2}^{*}$ and $r \equiv 1 \pmod{q+1}$.
Then~$f$ is a PP of $\F_{q^2}$ if and only if $(r, q^2-1)=1$ and $a^{q+1} \neq 1$.
If $f$ is a PP of $\F_{q^2}$ and $u$ is an inverse of\, $r$  modulo~$q^2-1$,
then its inverse on $\F_{q^2}$ is
\begin{equation*}
  f^{-1}(x) =  (1-a^{q+1})^{-u} (a^{q}x + x^{q})^{u}.
\end{equation*}
\end{corollary}

Next we give another example of the case $q^d -1 = s$ over $\F_{2^n}$.
\begin{corollary}\label{cor:qd-1=s}
Let $f(x)= x^r(x^{3} + a)^{5}$ where $a \in \F_{2^n}^{*}$, $n$ is even,
$r \equiv 1 \pmod{\ell}$ and $\ell = (2^n-1)/3$.
Then $f$ is a PP of $\F_{2^n}$ if and only if $(r, 2^n-1) = 1$ and $a^{\ell} \neq 1$.
If $f$ is a PP of $\F_{2^n}$ and $n \equiv 2 \pmod{4}$, then its inverse on $\F_{2^n}$ is
\[
f^{-1}(x) = a^{\ell u} x^{u}
\bigg( \sum_{i=1}^{n/2} a^{-\frac{16^i-1}{15}}x^{\frac{16^{i-1}-1}{5}}\bigg)^{5u},
\]
where $u$ is an inverse of\, $r$  modulo $2^{n}-1$.
\end{corollary}
\begin{proof}
  The permutation part is a direct consequence of \cref{thm:rst}.
  Let $\omega = a^{\ell}$. Then $\omega^3 = 1$ and, by $\omega \neq 1$,
  $\omega^2 + \omega + 1 = 0$. Thus $\omega/(1+\omega) = \omega^{-1}$
  and $\omega^{-5} = \omega$.
  Inserting them into~\cref{thm:qd-1-s-t} gives the above expression of $f^{-1}$.
\end{proof}

The next corollary is an example of the case $q^d -1 = t$ and $2n=3m$.
\begin{corollary}
Let $f(x)= x(x^{q+1} -a)^{q-1}$, where $a \in \F_{q^3}^{*}$ and $q$ is odd.
Then~$f$ is a PP of $\F_{q^3}$ if and only if $a^{(q^3-1)/2} = -1$.
In this case, its inverse on $\F_{q^3}$ is
\begin{equation*}
  f^{-1}(x) =  x \big(a^{q^2+q} + a^{q}x^{q+1} + x^{q^2+q+2}\big)^{q-1}.
\end{equation*}
\end{corollary}

\section{Permutation trinomials and tetranomials}\label{sec:PPs34}

Applying the main theorem to $t=2$, $3$, we can obtain
some permutation trinomials and tetranomials and their inverses,
which contain two classes of exceptional polynomials in \cref{thm:EPs}.
First, we give a simple lemma.

\begin{lemma}\label{gcd-lemma}
Let $a$ be odd, $m, n \in \mathbb{N}$ and $d=(m,n)$. Then
$(a^m-1)/(a^d -1)$ and $m/d$ have the same parity, and
\[
 \big( (a^m-1)/2, a^n-1 \big) =
\begin{cases}
    a^d -1     & \text{if $m/d$ is even,} \\
    (a^d -1)/2 & \text{if $m/d$ is odd.} \\
\end{cases}
\]
\end{lemma}
\begin{proof}
   Since $(a^m-1)/(a^d -1) =\sum_{i = 1}^{m/d}(a^d)^{i -1}$ and $a$ is odd,
   $(a^m-1)/(a^d -1)$ is an integer with the same parity as $m/d$. Then
  \[\begin{split}
    2\big( (a^m-1)/2, a^n-1 \big)
    &=  \big( a^m-1, 2(a^n-1) \big)  \\
    &= (a^m-1,a^n-1)\left(\frac{a^m-1}{(a^m-1,a^n-1)},
        \frac{2(a^n-1)}{(a^m-1,a^n-1)}\right) \\
    &= (a^m-1,a^n-1)\left(\frac{a^m-1}{(a^m-1,a^n-1)}, 2 \right) \\
    &= (a^d -1)\left(\frac{a^m-1}{a^d-1}, 2 \right) \\
    &= (a^d -1)(m/d, 2). \qedhere
  \end{split}\]
\end{proof}

Applying \cref{thm:rst} to $r=1$ and $t=2$, we derive the following result.

\begin{theorem}\label{q^m-t=2}
Let $f(x)= x^{q^m} -2 a x^{\frac{q^m+1}{2}} +a^2 x$,
where $a \in \F_{q^n}^{*}$, $q$ is odd and $m, n \in \mathbb{N}$.
Let $d = (m, n)$, $c = a^{(q^n-1)/(q^d-1)}$  and
\[
H_{2}(x) = x \bigg(\sum_{i=1}^{n/d} a^{-\frac{q^{i m}-1}{q^m-1}}
    x^{\frac{q^{(i-1) m}-1}{2}}\bigg)^{2}.
\]
\begin{enumerate}[\upshape(i)]
  \item If\, $m/d$ is even, then~$f$ is a PP of \,$\F_{q^n}$\! if and only if\, $c \neq 1$.
        In this case, 
        \[
        f^{-1}(x) = c^{2}(1-c)^{-2}\, H_{2}(x).
        \]
  \item If\, $m/d$ is odd, then~$f$ is a PP of \,$\F_{q^n}$\!
        if and only if\, $c^{2} \neq 1$.
        In this case, 
        \begin{equation}\label{f-1&r=1&t=2&odd}
        f^{-1}(x) = c^{2}(1-c^{2})^{-2}\big(2c x^{\frac{q^n-1}{2}} + c^{2} + 1\big) H_{2}(x).
        \end{equation}
\end{enumerate}
\end{theorem}
\begin{proof}
  Note that $f(x)=x(x^{\frac{q^m-1}{2}} - a)^2$.
  In the notation of \cref{thm:rst}, we have $t=2$ and $s = (q^m-1)/2$.
  If $m/d$ is even, then $(s, q^n-1) = q^d -1$ by \cref{gcd-lemma},
  and $\ell = (q^n-1)/(q^d-1)$. According to \cref{thm:rst},
  $f$ is a PP of $\F_{q^n}$ if and only if $c \neq 1$.
  It is easy to verify $f^{-1}(f(0)) = 0$.
  Next we only consider $x \in \F_{q^n}^{*}$ for computing $f^{-1}(x)$.
  Since $q^d-1 \mid s$, we have
  \begin{equation}\label{eq:N(xs)=1}
     (x^{s})^{\frac{q^n-1}{q^d-1}}
    = (x^{q^n-1})^{\frac{s}{q^d-1}}
    = 1
  \end{equation}
  for any $x \in \F_{q^n}^{*}$.
  Substituting \cref{eq:N(xs)=1} into \cref{eq:gx}, we obtain
  \begin{equation}\label{eq:gx+t=2+even}
    G(x) = c(1-c)^{-1}.
  \end{equation}

  If $m/d$ is odd, then $(s, q^n-1) = (q^d-1)/2$ by \cref{gcd-lemma},
  and $\ell = 2(q^n-1)/(q^d-1)$. According to \cref{thm:rst},
  $f$ is a PP of $\F_{q^n}$ if and only if $c^{2} \neq 1$.
  If $m/d$ is odd, then $(q^m-1)/(q^d-1)$ is also odd by \cref{gcd-lemma}, and so
  \[
      (x^{s})^{\frac{q^n-1}{q^d-1}}
    = x^{\frac{q^m-1}{2} \cdot \frac{q^n-1}{q^d-1}}
    = x^{\frac{q^n-1}{2} \cdot \frac{q^m-1}{q^d-1}}
    = x^{\frac{q^n-1}{2}}
  \]
  for any $x \in \F_{q^n}^{*}$.
  Since
  \[
      1 - c^2
     = \big(x^{\frac{q^n-1}{2}}    \big)^2 - c^2
     = \big(x^{\frac{q^n-1}{2}} + c\big)
       \big(x^{\frac{q^n-1}{2}} - c\big).
  \]
  We have
  \begin{equation}\label{eq:N(xs)2}
      \big(x^{\frac{q^n-1}{2}} - c\big)^{-1}
    =  \big(1 - c^{2}\big)^{-1}
      \big(x^{\frac{q^n-1}{2}} + c\big).
  \end{equation}
  Note that $c = a^{(q^n-1)/(q^d-1)}$.
  Inserting \cref{eq:N(xs)2} into \cref{eq:gx}, we obtain
  \begin{equation}\label{eq:gx+t=2+odd}
  G(x) =c(1 - c^{2})^{-1}
        (x^{\frac{q^n-1}{2}} + c).
  \end{equation}
  Then, for any $x \in \F_{q^n}^{*}$,
  \begin{equation}\label{eq:g1x2}
  \big(x^{\frac{q^n-1}{2}} + c\big)^{2}
  =  2c x^{\frac{q^n-1}{2}} + c^{2} + 1.
  \end{equation}
  Substituting \cref{eq:gx+t=2+even}, \cref{eq:gx+t=2+odd} and \cref{eq:g1x2}
  into \cref{thm:rst} gives the desire result.
\end{proof}

Taking $q^m = 5, 7, 9$ in \cref{q^m-t=2} leads to the following corollaries.

\begin{corollary}\label{cor:r=1+t=2+x5}
Let $f(x)=x^5 -2ax^3 +a^2x$, where $a \in \F_{5^n}^{*}$ and $n \in \mathbb{N}$.
Then $f$ is a PP of $\F_{5^n}$ if and only if $a^{(5^n-1)/2} =-1$.
In this case, the inverse of $f$ on $\F_{5^n}$ is
\[
f^{-1}(x)= 2a^{\frac{5^n-1}{4}}x^{\frac{5^n+1}{2}}
\bigg(\sum_{i=0}^{n-1} a^{-\frac{5^{i+1}-1}{4}}x^{\frac{5^{i}-1}{2}}\bigg)^{2}.
\]
\end{corollary}

\cref{cor:r=1+t=2+x5} is essentially \cite[Theorem 8]{Zheng7deg} or \cite[Lemma 4.9]{KLi18-1}.

\begin{corollary}\label{cor:r=1+t=2+x7}
Let $f(x)= x^{7} -2 a x^{4} +a^2 x$, where $a \in \F_{7^n}^{*}$ and $n \in \mathbb{N}$.
Then $f$ is a PP of $\F_{7^n}$ if and only if $a^{(7^n-1)/3} \neq 1$.
In this case, the inverse of $f$ on $\F_{7^n}$ is
\[
    f^{-1}(x)
    = 2x \big(2 a^{\frac{7^n-1}{6}} x^{\frac{7^n-1}{2}} + a^{\frac{7^n-1}{3}} + 1\big)
       \bigg(\sum_{i=0}^{n-1} a^{-\frac{7^{i+1}-1}{6}}x^{\frac{7^{i}-1}{2}}\bigg)^{2}.
\]
\end{corollary}
\begin{proof}
 The permutation part is a direct consequence of \cref{q^m-t=2}.
 Let $\omega = a^{(7^n-1)/3}$.
 Then $\omega^{3} = 1$ and $\omega \neq 1$ if $f$ is a PP.
 Hence $\omega^2 + \omega + 1 = 0$, and so $(1-\omega)^2 = -3\omega$.
 Inserting them into~\cref{f-1&r=1&t=2&odd} gives the above expression of $f^{-1}$.
\end{proof}

\begin{corollary}\label{cor:r=1+t=2+x9}
Let $f(x)=x^{9} + ax^{5} +a^2x$, where $a \in \F_{3^n}^{*}$ and $n \in \mathbb{N}$.
\begin{enumerate}[\upshape(i)]
  \item If $n$ is odd, then $f$ is a PP of $\F_{3^n}$ if and only if $a^{(3^n-1)/2} = -1$.
        In this case, the inverse of $f$ on $\F_{3^n}$ is
        \[
          f^{-1}(x)= x\bigg(\sum_{i=0}^{n-1} a^{-\frac{9^{i+1}-1}{8}}x^{\frac{9^{i}-1}{2}}\bigg)^{2}.
        \]
  \item If $n$ is even, then $f$ is a PP of $\F_{3^n}$ if and only if $c^{2} \neq 1$,
        where $c = a^{(3^n-1)/8}$. In this case, the inverse of $f$ on $\F_{3^n}$ is
        \begin{equation}\label{eq:f-1&x9&even}
          f^{-1}(x)= x \big(c^{5} x^{\frac{3^n-1}{2}} + c^{2} + 1\big)
          \bigg(\sum_{i=1}^{n/2} a^{-\frac{9^{i}-1}{8}}x^{\frac{9^{i-1}-1}{2}}\bigg)^{2}.
        \end{equation}
\end{enumerate}
\end{corollary}
\begin{proof}
 We only verify \cref{eq:f-1&x9&even}, because the rest parts
 are direct consequences of \cref{q^m-t=2}.
 Let $\omega = c^{2} = a^{(3^n-1)/4}$.
 Then $\omega^{4}=1$ and $\omega \neq 1$ if $f$ is a PP.
 Thus $\omega^{3} + \omega^{2} + \omega + 1 = 0$, and so
 $(1-\omega)^{2} = 1 + \omega + \omega^{2} = -\omega^{3}$ in $\F_{3^n}$.
 Then $\omega/(-\omega^{3}) = - \omega^{-2} = -\omega^{2}$
 and $-\omega^{2}(\omega + 1) = \omega + 1$.
 Inserting them into~\cref{f-1&r=1&t=2&odd} gives \cref{eq:f-1&x9&even}.
\end{proof}


Applying \cref{thm:rst} to $r=1$, $t=3$ and $q^m =7$
gives the following corollary.

\begin{corollary}\label{cor:r=1+t=3+x7}
Let $f(x)= x^{7} -3 a x^{5} +3a^2 x^3 -a^3 x$,
where $a \in \F_{7^n}^{*}$ and $n \in \mathbb{N}$.
Then $f$ is a PP of $\F_{7^n}$ if and only if $a^{(7^n-1)/2} = -1$.
In this case, its inverse on $\F_{7^n}$ is
\[
    f^{-1}(x)
    =x\big( 3 (a x^4)^{\frac{7^n-1}{6}}
      -3 (a x)^{\frac{7^n-1}{3}}  -2 \big)
      \bigg(\sum_{i=0}^{n-1} a^{-\frac{7^{i+1}-1}{6}}x^{\frac{7^{i}-1}{3}}\bigg)^{3}.
\]
\end{corollary}
\begin{proof}
  Clearly, $f(x)=x(x^{2} -a)^3$, and so $s=2$, $t=3$, $\ell = (7^n-1)/2$, $m = d = 1$.
  From \cref{thm:rst}, $f$ is a PP of $\F_{q^n}$ if and only if $a^{(7^n-1)/2} = -1$.
  It is easy to verify that $f^{-1}(f(0)) = 0$.
  Next we only consider $x \in \F_{q^n}^{*}$ for computing $f^{-1}(x)$.
  Note that
  \[
     2  = 1 - a^{\frac{7^n-1}{2}}
     = \big(  x^{\frac{7^n-1}{3}} \big)^{3}
     - \big(  a^{\frac{7^n-1}{6}} \big)^{3} \\
     =  \big(x^{\frac{7^n-1}{3}}
            - a^{\frac{7^n-1}{6}}\big)\lambda(x),
  \]
  where
  \[
     \lambda(x) = x^{\frac{2(7^n-1)}{3}}
          + a^{\frac{7^n-1}{6}}x^{\frac{7^n-1}{3}}
          + a^{\frac{7^n-1}{3}}.
  \]
  Therefore,
    \begin{equation}\label{eq:N(xs)3}
      \big(x^{\frac{7^n-1}{3}}
         - a^{\frac{7^n-1}{6}}\big)^{-1}
    =  4  \lambda(x).
  \end{equation}
  Inserting \cref{eq:N(xs)3} into \cref{eq:gx} yields
  $
    G(x)= 4 a^{\frac{7^n-1}{6}} \lambda(x).
  $
  Then, for any $x \in \F_{q^n}^{*}$,
\begin{equation}\label{eq:gx3}
  G(x)^3 = - \lambda(x)^3
  = 3(ax^4)^{\frac{7^n-1}{6}}
  - 3  (ax)^{\frac{7^n-1}{3}} -2.
\end{equation}
  Substituting \cref{eq:gx3} into \cref{thm:rst} gives the desire result.
\end{proof}

\section*{References}
\small
\bibliography{pps}

\begin{thebibliography}{40}
\providecommand{\natexlab}[1]{#1}
\providecommand{\url}[1]{\texttt{#1}}
\providecommand{\urlprefix}{URL }
\expandafter\ifx\csname urlstyle\endcsname\relax
  \providecommand{\doi}[1]{doi:\discretionary{}{}{}#1}\else
  \providecommand{\doi}[1]{doi:\discretionary{}{}{}\begingroup
  \urlstyle{rm}\url{#1}\endgroup}\fi
\providecommand{\bibinfo}[2]{#2}

\bibitem[{Akbary et~al.(2011)Akbary, Ghioca, and Wang}]{AGW}
\bibinfo{author}{A.~Akbary}, \bibinfo{author}{D.~Ghioca},
  \bibinfo{author}{Q.~Wang}, \bibinfo{title}{On constructing permutations of
  finite fields}, \bibinfo{journal}{Finite Fields Appl.} \bibinfo{volume}{17}
  (\bibinfo{year}{2011}) \bibinfo{pages}{51--67}.

\bibitem[{Akbary and Wang(2007)}]{AW07}
\bibinfo{author}{A.~Akbary}, \bibinfo{author}{Q.~Wang}, \bibinfo{title}{On
  polynomials of the form $x^rf(x^{(q-1)/l})$}, \bibinfo{journal}{Int. J. Math.
  Math. Sci.} \bibinfo{volume}{2007} (\bibinfo{year}{2007})
  \bibinfo{pages}{1--7, article ID 23408},
  \bibinfo{note}{\doi{10.1155/2007/23408}}.

\bibitem[{Bartoli(2020)}]{Bartoli20-tri}
\bibinfo{author}{D.~Bartoli}, \bibinfo{title}{Permutation trinomials over
  $\mathbb{F}_{q^3}$}, \bibinfo{journal}{Finite Fields Appl.}
  \bibinfo{volume}{61} (\bibinfo{year}{2020}) \bibinfo{pages}{101597}.

\bibitem[{Bartoli et~al.(2017)Bartoli, Giulietti, Quoos, and
  Zini}]{Bartoli8even}
\bibinfo{author}{D.~Bartoli}, \bibinfo{author}{M.~Giulietti},
  \bibinfo{author}{L.~Quoos}, \bibinfo{author}{G.~Zini},
  \bibinfo{title}{Complete permutation polynomials from exceptional
  polynomials}, \bibinfo{journal}{J. Number Theory} \bibinfo{volume}{176}
  (\bibinfo{year}{2017}) \bibinfo{pages}{46--66},
  \bibinfo{note}{\doi{10.1016/j.jnt.2016.12.016}}.

\bibitem[{Cao et~al.(2020)Cao, Hou, Mi, and Xu}]{CaoHMX20}
\bibinfo{author}{X.~Cao}, \bibinfo{author}{X.-D. Hou}, \bibinfo{author}{J.~Mi},
  \bibinfo{author}{S.~Xu}, \bibinfo{title}{More permutation polynomials with
  {Niho} exponents which permute $\mathbb{F}_{q^2}$}, \bibinfo{journal}{Finite
  Fields Appl.} \bibinfo{volume}{62} (\bibinfo{year}{2020})
  \bibinfo{pages}{101626}.

\bibitem[{Charpin et~al.(2016)Charpin, Mesnager, and Sarkar}]{Involutions}
\bibinfo{author}{P.~Charpin}, \bibinfo{author}{S.~Mesnager},
  \bibinfo{author}{S.~Sarkar}, \bibinfo{title}{Involutions over the {G}alois
  {F}ield $\mathbb{F}_{2^{n}}$}, \bibinfo{journal}{IEEE Trans. Inf. Theory}
  \bibinfo{volume}{62}~(\bibinfo{number}{4}) (\bibinfo{year}{2016})
  \bibinfo{pages}{2266--2276}.

\bibitem[{Coulter and Henderson(2004)}]{CoulterH04}
\bibinfo{author}{R.~Coulter}, \bibinfo{author}{M.~Henderson}, \bibinfo{title}{A
  note on the roots of trinomials over a finite field}, \bibinfo{journal}{Bull.
  Aust. Math. Soc.} \bibinfo{volume}{69} (\bibinfo{year}{2004})
  \bibinfo{pages}{429--432}.

\bibitem[{Coulter and Henderson(2002)}]{Coulter-1}
\bibinfo{author}{R.~S. Coulter}, \bibinfo{author}{M.~Henderson},
  \bibinfo{title}{The compositional inverse of a class of permutation
  polynomials over a finite field}, \bibinfo{journal}{Bull. Aust. Math. Soc.}
  \bibinfo{volume}{65} (\bibinfo{year}{2002}) \bibinfo{pages}{521--526}.

\bibitem[{Dickson(1896)}]{Dickson1896}
\bibinfo{author}{L.~E. Dickson}, \bibinfo{title}{The analytic representation of
  substitutions on a power of a prime number of letters with a discussion of
  the linear group, {P}art~{I}}, \bibinfo{journal}{Ann. Math.}
  \bibinfo{volume}{11} (\bibinfo{year}{1896}) \bibinfo{pages}{65--120},
  \bibinfo{note}{\url{https://www.jstor.org/stable/1967217}}.

\bibitem[{Fan(2019)}]{Fan7}
\bibinfo{author}{X.~Fan}, \bibinfo{title}{A classification of permutation
  polynomials of degree 7 over finite fields}, \bibinfo{journal}{Finite Fields
  Appl.} \bibinfo{volume}{59} (\bibinfo{year}{2019}) \bibinfo{pages}{1--21}.

\bibitem[{Fan(2020{\natexlab{a}})}]{Fan8even}
\bibinfo{author}{X.~Fan}, \bibinfo{title}{Permutation polynomials of degree 8
  over finite fields of characteristic 2}, \bibinfo{journal}{Finite Fields
  Appl.} \bibinfo{volume}{64} (\bibinfo{year}{2020}{\natexlab{a}})
  \bibinfo{pages}{101662}.

\bibitem[{Fan(2020{\natexlab{b}})}]{Fan8odd}
\bibinfo{author}{X.~Fan}, \bibinfo{title}{Permutation polynomials of degree 8
  over finite fields of odd characteristic}, \bibinfo{journal}{Bull. Aust.
  Math. Soc.} \bibinfo{volume}{101}~(\bibinfo{number}{1})
  (\bibinfo{year}{2020}{\natexlab{b}}) \bibinfo{pages}{40--55}.

\bibitem[{Hou et~al.(2020)Hou, Tu, and Zeng}]{HouTZ20-tri}
\bibinfo{author}{X.-D. Hou}, \bibinfo{author}{Z.~Tu},
  \bibinfo{author}{X.~Zeng}, \bibinfo{title}{Determination of a class of
  permutation trinomials in characteristic three}, \bibinfo{journal}{Finite
  Fields Appl.} \bibinfo{volume}{61} (\bibinfo{year}{2020})
  \bibinfo{pages}{101596}.

\bibitem[{Kyureghyan(2011)}]{Kyur11}
\bibinfo{author}{G.~M. Kyureghyan}, \bibinfo{title}{Constructing permutations
  of finite fields via linear translators}, \bibinfo{journal}{J. Combin. Theory
  Ser. A} \bibinfo{volume}{118} (\bibinfo{year}{2011})
  \bibinfo{pages}{1052--1061}.

\bibitem[{Li et~al.(2010)Li, Chandler, and Xiang}]{JLi10}
\bibinfo{author}{J.~Li}, \bibinfo{author}{D.~B. Chandler},
  \bibinfo{author}{Q.~Xiang}, \bibinfo{title}{Permutation polynomials of
  degree~6 or~7 over finite fields of characteristic~2},
  \bibinfo{journal}{Finite Fields Appl.} \bibinfo{volume}{16}
  (\bibinfo{year}{2010}) \bibinfo{pages}{406--419}.

\bibitem[{Li et~al.(2020)Li, Qu, Li, and Chen}]{KLiQLC20}
\bibinfo{author}{K.~Li}, \bibinfo{author}{L.~Qu}, \bibinfo{author}{C.~Li},
  \bibinfo{author}{H.~Chen}, \bibinfo{title}{On a conjecture about a class of
  permutation quadrinomials}, \bibinfo{journal}{Finite Fields Appl.}
  \bibinfo{volume}{66} (\bibinfo{year}{2020}) \bibinfo{pages}{101690}.

\bibitem[{Li et~al.(2018)Li, Qu, and Wang}]{KLi18}
\bibinfo{author}{K.~Li}, \bibinfo{author}{L.~Qu}, \bibinfo{author}{Q.~Wang},
  \bibinfo{title}{New constructions of permutation polynomials of the form $x^r
  h(x^{q-1})$ over $\mathbb{F}_{q^2}$}, \bibinfo{journal}{Des. Codes Cryptogr.}
  \bibinfo{volume}{86} (\bibinfo{year}{2018}) \bibinfo{pages}{2379--2405}.

\bibitem[{Li et~al.(2019)Li, Qu, and Wang}]{KLi18-1}
\bibinfo{author}{K.~Li}, \bibinfo{author}{L.~Qu}, \bibinfo{author}{Q.~Wang},
  \bibinfo{title}{Compositional inverses of permutation polynomials of the form
  $x^rh(x^s)$ over finite fields}, \bibinfo{journal}{Cryptogr. Commun.}
  \bibinfo{volume}{11} (\bibinfo{year}{2019}) \bibinfo{pages}{279--298}.

\bibitem[{Li and Helleseth(2019)}]{NLiH19}
\bibinfo{author}{N.~Li}, \bibinfo{author}{T.~Helleseth}, \bibinfo{title}{New
  permutation trinomials from {Niho} exponents over finite fields with even
  characteristic}, \bibinfo{journal}{Cryptogr. Commun.}
  \bibinfo{volume}{11}~(\bibinfo{number}{1}) (\bibinfo{year}{2019})
  \bibinfo{pages}{129--136}.

\bibitem[{Lidl and Niederreiter(1997)}]{FF}
\bibinfo{author}{R.~Lidl}, \bibinfo{author}{H.~Niederreiter},
  \bibinfo{title}{Finite Fields}, \bibinfo{publisher}{Cambridge {U}niv.
  {P}ress}, \bibinfo{year}{1997}.

\bibitem[{Muratovi\'{c}-Ribi\'{c}(2007)}]{MR-1}
\bibinfo{author}{A.~Muratovi\'{c}-Ribi\'{c}}, \bibinfo{title}{A note on the
  coefficients of inverse polynomials}, \bibinfo{journal}{Finite Fields Appl.}
  \bibinfo{volume}{13} (\bibinfo{year}{2007}) \bibinfo{pages}{977--980}.

\bibitem[{Niu et~al.(2020{\natexlab{b}})Niu, Li, Qu, and Wang}]{NiuLQW20}
\bibinfo{author}{T.~Niu}, \bibinfo{author}{K.~Li}, \bibinfo{author}{L.~Qu},
  \bibinfo{author}{Q.~Wang}, \bibinfo{title}{A general method for finding the
  compositional inverses of permutations from the AGW criterion},
  \bibinfo{howpublished}{arXiv:2004.12552},
  \bibinfo{note}{\url{https://arxiv.org/abs/2004.12552v1}},
  \bibinfo{year}{2020}{\natexlab{b}}.

\bibitem[{Niu et~al.(2020{\natexlab{a}})Niu, Li, Qu, and Wang}]{NiuLQW19}
\bibinfo{author}{T.~Niu}, \bibinfo{author}{K.~Li}, \bibinfo{author}{L.~Qu},
  \bibinfo{author}{Q.~Wang}, \bibinfo{title}{New constructions of involutions
  over finite fields}, \bibinfo{journal}{Cryptogr. Commun.}
  \bibinfo{volume}{12} (\bibinfo{year}{2020}{\natexlab{a}})
  \bibinfo{pages}{165--185}, \bibinfo{note}{\doi{10.1007/s12095-019-00386-2}}.

\bibitem[{Reis(2018)}]{Reis18}
\bibinfo{author}{L.~Reis}, \bibinfo{title}{Nilpotent linearized polynomials
  over finite fields and applications}, \bibinfo{journal}{Finite Fields Appl.}
  \bibinfo{volume}{50} (\bibinfo{year}{2018}) \bibinfo{pages}{279--292}.

\bibitem[{Shallue and Wanless(2013)}]{PPs6}
\bibinfo{author}{C.~J. Shallue}, \bibinfo{author}{I.~M. Wanless},
  \bibinfo{title}{Permutation polynomials and orthomorphism polynomials of
  degree six}, \bibinfo{journal}{Finite Fields Appl.} \bibinfo{volume}{20}
  (\bibinfo{year}{2013}) \bibinfo{pages}{84--92}.

\bibitem[{Tu et~al.(2019)Tu, Liu, and Zeng}]{TuLZ19}
\bibinfo{author}{Z.~Tu}, \bibinfo{author}{X.~Liu}, \bibinfo{author}{X.~Zeng},
  \bibinfo{title}{A revisit to a class of permutation quadrinomials},
  \bibinfo{journal}{Finite Fields Appl.} \bibinfo{volume}{59}
  (\bibinfo{year}{2019}) \bibinfo{pages}{57--85}.

\bibitem[{Tuxanidy and Wang(2014)}]{TW-1}
\bibinfo{author}{A.~Tuxanidy}, \bibinfo{author}{Q.~Wang}, \bibinfo{title}{On
  the inverses of some classes of permutations of finite fields},
  \bibinfo{journal}{Finite Fields Appl.} \bibinfo{volume}{28}
  (\bibinfo{year}{2014}) \bibinfo{pages}{244--281}.

\bibitem[{Tuxanidy and Wang(2017)}]{TW17}
\bibinfo{author}{A.~Tuxanidy}, \bibinfo{author}{Q.~Wang},
  \bibinfo{title}{Compositional inverses and complete mappings over finite
  fields}, \bibinfo{journal}{Discrete Appl. Math.} \bibinfo{volume}{217}
  (\bibinfo{year}{2017}) \bibinfo{pages}{318--329}.

\bibitem[{Wang(2009)}]{Wang-1}
\bibinfo{author}{Q.~Wang}, \bibinfo{title}{On inverse permutation polynomials},
  \bibinfo{journal}{Finite Fields Appl.} \bibinfo{volume}{15}
  (\bibinfo{year}{2009}) \bibinfo{pages}{207--213}.

\bibitem[{Wang(2017)}]{Wang-cyc2}
\bibinfo{author}{Q.~Wang}, \bibinfo{title}{A note on inverses of cyclotomic
  mapping permutation polynomials over finite fields}, \bibinfo{journal}{Finite
  Fields Appl.} \bibinfo{volume}{45} (\bibinfo{year}{2017})
  \bibinfo{pages}{422--427}.

\bibitem[{Wang(2019)}]{WangIndex19}
\bibinfo{author}{Q.~Wang}, \bibinfo{title}{Polynomials over finite fields: an
  index approach}, in: \bibinfo{editor}{K.-U. Schmidt},
  \bibinfo{editor}{A.~Winterhof} (Eds.), \bibinfo{booktitle}{Combinatorics and
  Finite Fields: Difference Sets, Polynomials, Pseudorandomness and
  Applications}, \bibinfo{pages}{319--346},
  \bibinfo{note}{\doi{10.1515/9783110642094-015}}, \bibinfo{year}{2019}.

\bibitem[{Wu(2013)}]{Wu-L-bi}
\bibinfo{author}{B.~Wu}, \bibinfo{title}{The compositional inverses of
  linearized permutation binomials over finite fields},
  \bibinfo{howpublished}{arXiv:1311.2154v1},
  \bibinfo{note}{\url{https://arxiv.org/abs/1311.2154}}, \bibinfo{year}{2013}.

\bibitem[{Wu and Liu(2013{\natexlab{b}})}]{Wu-bil}
\bibinfo{author}{B.~Wu}, \bibinfo{author}{Z.~Liu}, \bibinfo{title}{The
  compositional inverse of a class of bilinear permutation polynomials over
  finite fields of characteristic $2$}, \bibinfo{journal}{Finite Fields Appl.}
  \bibinfo{volume}{24} (\bibinfo{year}{2013}{\natexlab{b}})
  \bibinfo{pages}{136--147}.

\bibitem[{Wu and Liu(2013{\natexlab{a}})}]{Wu-L}
\bibinfo{author}{B.~Wu}, \bibinfo{author}{Z.~Liu}, \bibinfo{title}{Linearized
  polynomials over finite fields revisited}, \bibinfo{journal}{Finite Fields
  Appl.} \bibinfo{volume}{22} (\bibinfo{year}{2013}{\natexlab{a}})
  \bibinfo{pages}{79--100}.

\bibitem[{Zha et~al.(2019)Zha, Hu, and Zhang}]{ZhaHZ19}
\bibinfo{author}{Z.~Zha}, \bibinfo{author}{L.~Hu}, \bibinfo{author}{Z.~Zhang},
  \bibinfo{title}{Permutation polynomials of the form $x + \gamma
  \mathrm{Tr}_{q}^{q^n}(h(x))$}, \bibinfo{journal}{Finite Fields Appl.}
  \bibinfo{volume}{60} (\bibinfo{year}{2019}) \bibinfo{pages}{101573}.

\bibitem[{Zheng et~al.(2019)Zheng, Yuan, Li, Hu, and Zeng}]{ZhengYL+2019}
\bibinfo{author}{D.~Zheng}, \bibinfo{author}{M.~Yuan}, \bibinfo{author}{N.~Li},
  \bibinfo{author}{L.~Hu}, \bibinfo{author}{X.~Zeng},
  \bibinfo{title}{Constructions of involutions over finite fields},
  \bibinfo{journal}{IEEE Trans. Inf. Theory}
  \bibinfo{volume}{65}~(\bibinfo{number}{12}) (\bibinfo{year}{2019})
  \bibinfo{pages}{7876--7883}, \bibinfo{note}{\doi{10.1109/TIT.2019.2919511}}.

\bibitem[{Zheng et~al.(2020{\natexlab{b}})Zheng, Wang, Wang, and
  Wei}]{ZhengRDP}
\bibinfo{author}{Y.~Zheng}, \bibinfo{author}{F.~Wang},
  \bibinfo{author}{L.~Wang}, \bibinfo{author}{W.~Wei}, \bibinfo{title}{On
  inverses of some permutation polynomials over finite fields of characteristic
  three}, \bibinfo{journal}{Finite Fields Appl.} \bibinfo{volume}{66}
  (\bibinfo{year}{2020}{\natexlab{b}}) \bibinfo{pages}{101670},
  \bibinfo{note}{\doi{10.1016/j.ffa.2020.101670}}.

\bibitem[{Zheng et~al.(2020{\natexlab{a}})Zheng, Wang, and Wei}]{Zheng7deg}
\bibinfo{author}{Y.~Zheng}, \bibinfo{author}{Q.~Wang},
  \bibinfo{author}{W.~Wei}, \bibinfo{title}{On inverses of permutation
  polynomials of small degree over finite fields}, \bibinfo{journal}{IEEE
  Trans. Inf. Theory} \bibinfo{volume}{66}~(\bibinfo{number}{2})
  (\bibinfo{year}{2020}{\natexlab{a}}) \bibinfo{pages}{914--922},
  \bibinfo{note}{\doi{10.1109/TIT.2019.2939113}}.

\bibitem[{Zheng et~al.(2016)Zheng, Yu, Zhang, and Pei}]{ZhengPW2}
\bibinfo{author}{Y.~Zheng}, \bibinfo{author}{Y.~Yu},
  \bibinfo{author}{Y.~Zhang}, \bibinfo{author}{D.~Pei},
  \bibinfo{title}{Piecewise constructions of inverses of cyclotomic mapping
  permutation polynomials}, \bibinfo{journal}{Finite Fields Appl.}
  \bibinfo{volume}{40} (\bibinfo{year}{2016}) \bibinfo{pages}{1--9}.

\bibitem[{Zheng et~al.(2015)Zheng, Yuan, and Pei}]{ZhengPW1}
\bibinfo{author}{Y.~Zheng}, \bibinfo{author}{P.~Yuan},
  \bibinfo{author}{D.~Pei}, \bibinfo{title}{Piecewise constructions of inverses
  of some permutation polynomials}, \bibinfo{journal}{Finite Fields Appl.}
  \bibinfo{volume}{36} (\bibinfo{year}{2015}) \bibinfo{pages}{151--169}.

\end{thebibliography}

\end{document}